\documentclass[12pt,english]{amsart}

\input{Dynkin.sty}
\usepackage{amsmath,amsfonts,amssymb,amsthm,amscd}
\usepackage{dsfont, mathrsfs,enumerate,eucal,fontenc}
\usepackage{xspace}

\usepackage[english]{babel}

\usepackage{xy, xypic}
\usepackage{color} 
\usepackage{graphicx}
\usepackage{psfrag}  
\usepackage{pst-all}

\usepackage{youngtab}

\usepackage{xcolor}
\definecolor{rouge}{rgb}{0.8,0.1,.4}
\definecolor{bleu}{rgb}{0.1,0.2,0.9}
\definecolor{violet}{rgb}{0.7,0,0.8}
\usepackage[colorlinks=true,linkcolor=bleu,urlcolor=violet,citecolor=rouge]{hyperref}
\usepackage{breakurl}
\usepackage{xspace}

\newtheorem{thm}{Theorem}[section]
\newtheorem{lemma}[thm]{Lemma}

\theoremstyle{definition}
\newtheorem{defi}[thm]{Definition}

\newtheorem{ex}[thm]{Example}

\usepackage{fullpage}

\def\le{\leqslant}
\def\ge{\geqslant}

 \def\g{{\mathfrak{g}}}  % ground alg de Lie
\def\l{{\mathfrak{l}}}  % Levi

\def\N{\mathbb N}

\def\eps{\varepsilon}

\def\P{\mathcal P}
\def\O{\mathcal  O}      % \C [[ t ]]
\def\Ss{\mathcal S}

\def\PT{{\rm PT}}
\def\Mo{{\rm M}}
\def\i{{\bf i}}
\def\zPT{z_{\PT}}
\def\zM{z_{\Mo}}

\begin{document}

\title[Corrigendum to "On the Dimension of the Sheets of a Reductive Lie Algebra"]
{Corrigendum to "On the Dimension of the Sheets of a Reductive Lie Algebra"}

\author[Anne Moreau]{Anne Moreau}

\address{\SMALL Anne Moreau, Laboratoire de Math{\'e}matiques et Applications, 
 Universit{\'e} de Poitiers, France}
 
\email{anne.moreau@math.univ-poitiers.fr}

\subjclass{14A10, 14L17, 22E20, 22E46}

\keywords
{coadjoint orbit, sheet, induced nilpotent orbit, rigid nilpotent orbit}

\begin{abstract}
This note is a corrigendum to \cite{M08}. 
As it has been recently pointed out to me by Alexander Premet, 
\cite[Remark 3.12]{M08} is incorrect. 
%An error in \cite[Proposition 3.11]{M08} was recently pointed out to me by Alexander Premet. 
We explain in this note the impacts of that error in \cite{M08}, and amend 
certain of its statements. 
In particular, we verify that the statement of \cite[Theorem 3.13]{M08} remains correct 
in spite of this error. 
\end{abstract}

\date{\today}

\maketitle

\section{Introduction}
Let $\g$ be a complex simple Lie algebra and $G$ its adjoint group. 
We investigate in \cite{M08} the dimension of the subsets, for $m \in \N$, 
$$\g^{(m)}:=\{x \in \g \; | \; \dim (Gx)=2m\},  
$$
where $Gx$ denotes the adjoint orbit of $x\in \g$. 
The irreducible components of the subsets $\g^{(m)}$ 
are called the {\em sheets} of $\g$, \cite{BoK,Bor}. 
Thus, for any $m \in \N$, 
\begin{eqnarray}   \label{eq:gm}
\dim \g^{(m)} = \max \{\dim \Ss \; ; \; \Ss \subset \g^{(m)}\} ,
\end{eqnarray}
where $\Ss$ runs through all sheets contained 
in $\g^{(m)}$. 
The sheets are known to be
parameterized by the pairs $(\l,\O_{\l})$,  
up to $G$-conjugacy class, consisting of a Levi subalgebra
$\l$ of $\g$ and a rigid nilpotent orbit $\O_\l$ in $\l$, cf.\ \cite{Bor}. 
This parametrization enables to write the dimension 
of a sheet $\Ss$ associated with a pair $(\l,\O_{\l})$ 
as the sum of the dimension of the center of $\l$  
and the dimension of the unique 
nilpotent orbit contained in $\Ss$, see e.g.\ \cite[Proposition 2.11]{M08}. 

In the classical case, formulas for $\g^{(m)}$ 
are given in  \cite[Theorems 3.3 and 3.13]{M08} in term of partitions associated with nilpotent 
elements of $\g$. 
As it has been recently pointed out 
by Alexander Premet, %Proposition 3.11 and 
Remark 3.12 in \cite{M08} 
which claims that  {\em "in the classical case, 
the dimension of a sheet containing a given nilpotent 
orbit does not depend on the choice of a sheet containing it"} 
is incorrect.  
We give here some counter-examples  
(cf.\ Examples \ref{ex1} and \ref{ex2}; see also \cite[Remark 4]{PT}).  
This is true only for the type {\bf A} where each nilpotent 
element belongs to only one sheet. 
The error stems from the proof of \cite[Proposition 3.11]{M08}; 
see Section \ref{S:err} for explanations. 
As a consequence, the proof of \cite[Theorems 3.13]{M08},  
partly based on \cite[Proposition 3.11]{M08}, is incorrect too.   
However its statement remains true. 
This can be shown through a recent work of Premet and Topley, 
\cite{PT}. 
In more details, another formula for $\g^{(m)}$ 
in term of partitions can be traced out from \cite[Corollary 9]{PT} 
and the equality (\ref{eq:gm}). 
In this note, we verify (cf.\ Theorems \ref{t:z})  
that the Premet-Topley formula for $\g^{(m)}$ coincides 
with the one of \cite[Theorem 3.13]{M08}.  

\smallskip

The note is organized as follows. 

\smallskip

In Section \ref{S:pf}, we recall some definitions and results 
of \cite{PT} and show that the statement of \cite[Theorem 3.13]{M08} is correct 
in spite of the error in \cite[Proposition 3.11]{M08}, see Theorem \ref{t:z}(ii). 
In Section \ref{S:err}, we precisely pin down  the error 
in the proof \cite[Proposition 3.11]{M08} and 
describe the impacts of that error in \cite{M08}. 
As a conclusion, we list in Section \ref{S:concl} 
all corrections which have to be taken into account in \cite{M08}. 

\smallskip 

Since the corrections in \cite{M08} only concern 
the types {\bf B}, {\bf C} and {\bf D}, we assume 
for the remaining of the note that 
$\g$ is either $\mathfrak{so}(N)$ 
or $\mathfrak{sp}(N)$, with $N \ge 2$, and $\eps$ is 1 or $-1$ 
depending on whether $\g =\mathfrak{so}(N)$ or $\mathfrak{sp}(N)$.  
Following the notations of \cite{M08} (or \cite{PT}), 
we denote by $\P_\eps(N)$ the set of partitions of $N$ associated 
with the nilpotent elements of $\g$. 
For $\lambda=(\lambda_1, \ldots,\lambda_n) \in \P_\eps(N)$,  
we denote by $e(\lambda)$ the corresponding nilpotent element of $\g$ 
whose Jordan block sizes are $\lambda_1, \ldots,\lambda_n$. 
We will always assume that $\lambda_1 \ge \cdots \ge \lambda_n$. 

\bigskip

\noindent
{\bf Acknowledgments.} I would like to thank A. Premet for having pointed 
out to me the error in my paper, and Lewis Topley for useful discussions 
and explanations.  
I also take the opportunity to thank Oscar Chacaltana for his 
interest in the subject and interesting e-mail exchanges.

\section{The main result}  \label{S:pf}
For the convenience of the reader, we recall here all the necessary definitions 
and results of \cite{PT}. 
Given a partition $\lambda =(\lambda_1,\ldots,\lambda_n) \in \P_\eps(N)$ 
we set, 
$$\Delta(\lambda) := \{ 1 \le i < n \; ; \; \eps (-1)^{\lambda_1}  = \eps (-1)^{\lambda_{i+1}} = -1, 
\; \lambda_{i-1} \not= \lambda_i \ge \lambda_{i+1} \not= \lambda_{i+2}\}. 
$$
Our convention is that $\lambda_0 =0$ and $\lambda_i =0$ for all $i > n$. 
Recall the following result of Kempken and Spaltenstein 
(also recalled in \cite{M08} and \cite{PT}):

\begin{thm}[{\cite{Ke,Sp}}]  \label{t:rig}
Let $\lambda=(\lambda_1,\ldots, \lambda_n) \in \P_\eps(N)$. 
Then $e(\lambda)$ is rigid if and only if 
\begin{itemize}
\item $\lambda_i - \lambda_{i+1} \in \{0,1\}$ for all $1 \le i \le n$;
\item the set $\{i \in \Delta(\lambda) \; ; \; \lambda_i = \lambda_{i+1} \}$ 
is empty. 
\end{itemize} 
\end{thm}

Denote by $\P_\eps^*(N)$ the set of $\lambda \in \P_\eps(N)$ 
such that $e(\lambda)$ is rigid.  
We call the elements of $\P_\eps^*(N)$ the {\em rigid partitions}. 
We first introduce the notion of {\em admissible sequences}, see \cite[\S3.1]{PT}. 
This is an extended version of the algorithm described in \cite{M08} 
which takes $\lambda\in\P_\eps(N)$ and returns an element of $\P_\eps^*(N)$ 
compatible for the induction process of nilpotent orbits. 

\smallskip

Let $\i$ be a finite sequence of integers between $1$ and
$n$. 
The procedure of \cite{PT} is as follows: 
the algorithm commences with input $\lambda = \lambda^\i \in \P_\eps(N)$ 
where $\i = \varnothing$ is the empty sequence. 
At the $l^{\rm th}$ iteration, the algorithm takes 
$\lambda^\i  \in \P_\eps(N - 2 \sum_{j=1}^{l-1} i_j)$ 
where $\i =(i_1,\ldots,i_{l-1})$ and returns $\lambda^{\i'} \in  \P_\eps(N - 2 \sum_{j=1}^{l} i_j)$ 
where $\i' = (i_1,\ldots,i_{l-1},i_l)$ for some $i_l$. 
If the output $\lambda^{\i'}$ is a rigid partition then the algorithm 
terminates after the $l^{\rm th}$ iteration with output $\lambda^{\i'}$ . 
We shall now explicitly describe the $l^{\rm th}$ iteration of the algorithm. 
If after the $(l-1)^{\rm th}$ iteration the input $\lambda^{\i}$ is not rigid 
then the algorithm behaves as follows. 
Let $i_l$ denote any index in the range $1 \le i \le n$ 
such that either of the following case occur:

\begin{description}
\item[Case 1]  $\quad \lambda_{i_l}^{\i} \ge \lambda_{i_l+1}^{\i} +2$;
\smallskip
\item[Case 2]   $\quad i_l \in \Delta(\lambda^{\i})$ and $\lambda_{i_l}^{\i}
 = \lambda_{i_l+1}^{\i}$.
\end{description}

\noindent
Note that no integer $i_l$ will fulfill both criteria. 
If $\i=(i_1,\ldots,i_{l-1})$ then define $\i' =(i_1,\ldots,i_{l-1},i_l)$. 
For Case 1 the algorithm has output 
$$\lambda^{\i'} = (\lambda_1^{\i}-2,\lambda_2^{\i}-2,\ldots, \lambda_{i_l}^{\i}-2, 
\lambda_{i_l+1}^{\i},\ldots, \lambda_{n}^{\i} )
$$
whilst for Case 2 the algorithm has output 
$$\lambda^{\i'} = (\lambda_1^{\i}-2,\lambda_2^{\i}-2,\ldots, \lambda_{i_l-1}^{\i}-2, 
\lambda_{i_l}^{\i}-1, \lambda_{i_l+1}^{\i}-1, 
\lambda_{i_l+2}^{\i},\ldots, \lambda_{n}^{\i} ). 
$$
Due to its definition and the classification of rigid partitions the above algorithm 
certainly terminates after a finite number of steps.

\begin{defi}[{\cite[\S3.1]{PT}}]   \label{d:adm} 
We say that a sequence $\i =(i_1,\ldots,i_l)$ is an {\em admissible sequence} 
for $\lambda$ if Case 1 or Case 2 occurs at the point $i_k$ for the partition
$\lambda^{(i_1,\ldots,i_{k-1})}$ for each $k=1,\ldots,l$. 
An admissible sequence $\i$ for $\lambda$ is be called {\em a maximal admissible sequence
 for $\lambda$} if neither Case 1 nor Case 2 occurs for any index $i$ between $1$ and $n$ 
 for the partition $\lambda^\i$. 
 By convention the empty sequence is admissible for any $\lambda \in \P_\eps(N)$.
 \end{defi}

\noindent
As observed in \cite[Lemma 6]{PT}, if $\i$ is an admissible sequence for $\lambda$,  
then $\i$ is maximal admissible if and only if $\lambda^{\i}$ is a rigid partition. 
We will denote by $|\i| :=l$ the length of an admissible sequence for $\lambda$.

\begin{defi} \label{d:can}
The algorithm as described in \cite{M08} corresponds to the special 
case where in the above algorithm, we define at each step   
$i_l$ to be the smallest integer which fulfills 
one the Case 1 or Case 2 criteria, and $\lambda^\i$ is rigid. 
In the sequel, we will refer the so obtained maximal admissible 
sequence for $\lambda$ to as the {\em canonical maximal admissible 
sequence for $\lambda$} and we denote it by $\i^0$.  
Then we set 
$$z_\Mo (\lambda):= |\i^0|.$$
\end{defi}

\noindent
{\em Remark.} The integer $\zM(\lambda)$ corresponds 
to the integer $z(\lambda)$ of \cite{M08}.

\begin{defi}[{\cite[Definition 1]{PT}}]  \label{d:2step} 
If $i \in \Delta(\lambda)$ then the pair $(i,i+1)$ is called a {\em 2-step of $\lambda$}. 
If $i > 1$ and $(i,i+1)$ is a 2-step of $\lambda$ then $\lambda_{i-1}$ and 
$\lambda_{i+2}$ are referred to as the {\em boundary of $(i,i+1)$}. 
If $1 \in \Delta(\lambda)$ then $\lambda_3$ is referred to as the boundary of $(1,2)$ 
(if $n = 2$ then $\lambda_3 = 0$ by convention). 
\end{defi}

We observe that $\Delta(\lambda)$ is the set of 2-steps of $\lambda$, 
and by $|\Delta(\lambda)|$ its cardinality.

\begin{defi}[{\cite[\S3.2]{PT}}]  \label{d:bad2step}
If $i \in \Delta(\lambda)$ then we say that the 2-step $(i, i + 1)$ has a {\em good boundary} 
if $\lambda_1$ and the boundary of $(i, i + 1)$ have the opposite parity. 
If the boundary of a 2-step $(i, i + 1)$ of $\lambda$ is not good then we say that it is {\em bad} 
and we refer to $(i,i+1)$ as a {\em bad 2-step}. 
Note that $(i,i+1)$ is a bad 2-step of $\lambda$ 
if and only if either $i>1$ and $\lambda_{i-1} - \lambda_i \in 2\N$, 
or $\lambda_{i+1} - \lambda_{i+2} \in 2\N$.
\end{defi}

We denote by $\Delta_{\rm bad}(\lambda)$ the set 
of bad 2-steps of $\lambda$, and by $|\Delta_{\rm bad}(\lambda)|$ 
its cardinality.

\begin{defi}  [{\cite[Definition 2]{PT}}]   \label{d:cluster}
A sequence $1 \le i_1 < \cdots  < i_k < n$ with 
$k \ge 2$ is called a {\em 2-cluster of $\lambda$} 
whenever $i_j \in \Delta(\lambda)$ and $i_{j+1} = i_j +2$ for all $j$. 
We say that a 2-cluster $i_1,\ldots, i_k$ has a {\em bad boundary} 
if either of the following conditions holds:
\begin{itemize}
\item $\lambda_{i_1-1} - \lambda_{i_1} \in 2\N$;
\item $\lambda_{i_k+1} - \lambda_{i_k+2} \in 2\N$. 
\end{itemize}
(if $i_1 = 1$ the the first condition should be omitted). 
A {\em bad 2-cluster} is one which has a bad boundary, 
whilst a {\em good 2-cluster} is one without a bad boundary. 
\end{defi}

We denote by  $\Sigma(\lambda)$ the set of good 2-clusters of $\lambda$, 
and by $|\Sigma(\lambda)|$ its cardinality.

\begin{lemma}[{\cite[Lemma 11]{PT}}]  \label{l:cluster}
A good 2-cluster is maximal in the sense that it is not a proper 
subsequence of any 2-cluster. 
\end{lemma}

\begin{defi}[Premet-Topley] \label{d:zPT} 
For any $\lambda \in \P_\eps(\lambda)$, the integer 
$\zPT(\lambda)$ is defined by the formula: 
$$\zPT(\lambda ) := s(\lambda) + |\Delta(\lambda)| - 
	|\Delta_{\rm bad}(\lambda)| + |\Sigma(\lambda)|$$
where 
$$s(\lambda) := \sum_{i=1}^{n} 
\left[ (\lambda_i - \lambda_{i+1})/2\right].$$
\end{defi}

\noindent
{\em Remark.} The integer $\zPT(\lambda)$ corresponds 
to the integer $z(\lambda)$ of \cite{PT}. 

\medskip

\noindent
By \cite[Theorem 8]{PT}, we have that 
\begin{eqnarray}  \label{eq:zPT}
\zPT(\lambda) := \max |\i|
\end{eqnarray}
where the maximum is taken over all admissible sequences for $\lambda$. 
Hence, by \cite[Corollary 9]{PT} and the equality (\ref{eq:gm}) 
of the introduction, we get:

\begin{thm}[Premet-Topley]  \label{t:PT}
For any $m \in \N$, we have
$$\dim \g^{(m)} = 2m + \max \{ \zPT(\lambda) \; ; \; \lambda \in \P_\eps(N)
\textrm{ s.t } \dim G e(\lambda) = 2m\}.
$$
\end{thm}

The main result of this note is:

\begin{thm}  \label{t:z}
{\rm (i)} 
For any $\lambda \in \P_\eps(N)$, 
we have  $z_\Mo(\lambda)= \zPT(\lambda)$. 

\noindent
{\rm (ii)} For any $m \in \N$, we have
$$\dim \g^{(m)} = 2m + \max \{ z_\Mo (\lambda) \; ; \; \lambda \in \P_\eps(N) 
\textrm{ s.t } \dim G e(\lambda) = 2m\}.
$$
In other words, the statement of \cite[Theorem 3.13]{M08} is correct. 
\end{thm}

\begin{proof}
(ii) is a direct consequence of (i) and Theorem \ref{t:PT}. 

\smallskip

\noindent
(i) We argue by induction on $N$ 
(the statement is true for small $N$). 
Let $N > 2$ and assume the statement true for any 
$\lambda \in \P_\eps(N')$, with $1 \le N' \le N$,  
and let $\lambda \in \P_\eps(N)$. 

If $\lambda \in \P_\eps^*(N)$, then $\zPT(\lambda) = z_\Mo(\lambda) =0$ 
(see Theorem \ref{t:rig}, Definition \ref{d:adm} and equality (\ref{eq:zPT})). 
So, we can assume that $\lambda$ is not a rigid partition. 
In particular, $\zPT(\lambda) >0$ and $z_\Mo(\lambda) >0$. 
To ease notation, we simply denote here by $\i := \i^0$ the canonical 
maximal sequence for $\lambda$. 
Then recall that by Definition \ref{d:can}, $z_\Mo(\lambda) = |\i|$.  
Set $\lambda' := \lambda^{(i_1)}$. 
Clearly, $z_\Mo(\lambda') = z_\Mo(\lambda) -1$. 
By the induction hypothesis, we have $\zPT(\lambda') = z_\Mo(\lambda')$. 
Hence, we have to show that: 
$$\zPT(\lambda') = \zPT(\lambda) -1 .
$$ 
Our strategy is to compare the formulas for $\zPT(\lambda')$ and $\zPT(\lambda)$ 
given by Definition \ref{d:zPT}. 
Recall that $i_1$ is the smallest integer which 
fulfills one of the Case 1 or Case 2 criteria for $\lambda$. 
First of all, we observe that if $i \in \Delta(\lambda)$ (resp.\ $i \in \Delta(\lambda')$),  
then $i \ge i_1$.  
Indeed, if $i \in \Delta(\lambda)$ and $i < i_1$ (if $i_1=1$, it is clear), 
then either $\lambda_i = \lambda_{i+1}$ and then $i$ fulfills the 
Case 2 which contradicts the minimality of $i_1$, 
or $\lambda_i - \lambda_{i+1} \in 2\N \smallsetminus\{0\}$ 
and then $i$ fulfills the 
Case 1 which contradicts the minimality of $i_1$ too. 

We now consider the two situations Case 1 and Case 2 separately. 

\smallskip

\noindent
\underline{\bf Case 1}: {$\lambda_{i_1} \ge \lambda_{i_1 +1} +2$}.

\smallskip

\noindent
We have, 
$$\lambda' = (\lambda_1 -2, \ldots, \lambda_{i_1-1}-2, \lambda_{i_1} -2, 
\lambda_{i_1+1},\ldots, \lambda_n), 
$$
and 
\begin{eqnarray*}
s(\lambda') &=& \displaystyle{\sum_{i=1}^{i_1-1} 
\left[ (\lambda_i  - \lambda_{i+1} ) /2 \right]
+\left[  (\lambda_{i_1} - 2 - \lambda_{i_1+1})/2\right]
+ \sum_{i=i_1+1}^{n} 
\left[ ( \lambda_i  - \lambda_{i+1})/2 \right]
} \\
&=& s(\lambda)-1.
\end{eqnarray*}
Compare now the other terms appearing in Definition \ref{d:zPT}. 
Note that $i_1 \in \Delta(\lambda)$ (resp. 
$i_1\in \Delta_{\rm bad}(\lambda)$)
if and only if $i_1 \in \Delta(\lambda')$ (resp. $i_1\in \Delta_{\rm bad}(\lambda')$)
since the passing from $\lambda$ to $\lambda'$ preserves 
the parities. 
For the same reason, $i_1$ belongs to a good 2-cluster 
of $\lambda$ if and only $i_1$ belongs to a good 2-cluster 
of $\lambda'$. 
%Also, if $i_1+1 \in \Delta(\lambda')$ (resp. $\Delta_{\rm bad}(\lambda')$), 
%then $i_1+1 \in \Delta(\lambda)$ (resp. $\Delta_{\rm bad}(\lambda)$). 

Then we discuss two cases depending on whether 
$i_1+1$ is in $\Delta(\lambda)$ or not:

\smallskip

\noindent
\textbullet \;  {$i_1+1 \in \Delta(\lambda)$}.\\
Once again, we consider two cases:

\begin{itemize}
\item[\textasteriskcentered] {$\lambda_{i_1}-2 \not=\lambda_{i_1+1}$}.\\  
Then $i_1+1 \in \Delta(\lambda')$ too. 
Moreover, $i_1+1 \in \Delta_{\rm bad}(\lambda')$ 
if and only if $i_1+1 \in \Delta_{\rm bad}(\lambda)$. 
Hence, we conclude that $|\Delta(\lambda')| = |\Delta(\lambda)|$, 
$|\Delta_{\rm bad}(\lambda')| = |\Delta_{\rm bad}(\lambda)|$  and 
$|\Sigma(\lambda')|=|\Sigma(\lambda)|$. 

\item[\textasteriskcentered]  {$\lambda_{i_1}-2 =\lambda_{i_1+1}$}.\\
Then $i_1+1 \in \Delta_{\rm bad}(\lambda)$ 
since $\lambda_{i_1} - \lambda_{i_1+1} =2 \in 2\N$. 
But $i_1+1 \not\in \Delta(\lambda')$. 
Therefore, $|\Delta(\lambda')| = |\Delta(\lambda)| -1$ 
and $|\Delta_{\rm bad}(\lambda')| = |\Delta_{\rm bad}(\lambda)|-1$. 
Moreover, if $i_1+1$ belongs to a 2-cluster of $\lambda$, 
then it is bad because $\lambda_{i_1} - \lambda_{i_1+1} \in 2\N$. 
Hence, we have $|\Sigma(\lambda')|=|\Sigma(\lambda)|$.
\end{itemize}

\smallskip

\noindent
\textbullet \;  {$i_1+1 \not\in \Delta(\lambda)$}.\\
In this case, note that $i_1+1 \not\in \Delta(\lambda')$. 
Hence, we conclude that $|\Delta(\lambda')| = |\Delta(\lambda)|$, 
$|\Delta_{\rm bad}(\lambda')| = |\Delta_{\rm bad}(\lambda)|$  and 
$|\Sigma(\lambda')|=|\Sigma(\lambda)|$. 

\medskip

\noindent
\underline{\bf Case 2}: {$i_1 \in \Delta(\lambda)$ and 
$\lambda_{i_1} = \lambda_{i_1+1}$}.

\smallskip

\noindent
By the minimality condition of $i_1$, we have $\lambda_{i_1-1} = \lambda_{i_1}+1$ 
(except for $i_1=1$, 
in which case $\lambda_{i_1-1} =0$ by convention), 
and so $\lambda_{i_1-2} = \lambda_{i_1-1}$ 
because $\eps(-1)^{ \lambda_{i_1-1}} = 1$.  
We have 
$$\lambda' = (\lambda_1 -2, \ldots, \lambda_{i_1-1}-2, \lambda_{i_1} -1, 
\lambda_{i_1+1} -1, \lambda_{i_1 +2},\ldots, \lambda_n), 
$$
and
\begin{eqnarray*}
s(\lambda') &=& \sum_{i=1}^{i_1 -2} 
\left[ (\lambda_i - \lambda_{i+1})/2 \right]
+ \underbrace{
\left[(\lambda_{i_1-1} - \lambda_{i_1} -1 )/2 \right]
}
_{=0 %=\left[(\lambda_{i_1-1} - \lambda_{i_1} )/2 \right] \atop 
\textrm{ since }\lambda_{i_1-1} = \lambda_{i_1}+1} \\
&& + \ \left[ (\lambda_{i_1}  - \lambda_{i_1+1})/2 \right] + 
\left[ \lambda_{i_1+1} - 1 - \lambda_{i_1+2})/2 \right] 
+ \sum_{i=i_1+1}^{n} 
\left[ (\lambda_i - \lambda_{i+1} ) /2 \right] \\
&=& \begin{cases} 
s(\lambda)-1& \textrm{ if } \lambda_{i_1+1} - \lambda_{i_1+2} \in 2\N ;\\
s(\lambda)& \textrm{ if } \lambda_{i_1+1}- \lambda_{i_1+2} \not\in 2\N.   
\end{cases}
\end{eqnarray*}
(If $i_1=0$, we start at the second line and we get the same conclusion.) 
Also, observe that in Case 2, we have 
$$|\Delta(\lambda')| = |\Delta(\lambda)|-1 .$$ 
Indeed, $i_1 \in \Delta(\lambda)$ but $i_1 \not\in \Delta(\lambda')$ 
and for the indexes $i \not =i_1$ we have here the equivalence: 
$i \in \Delta(\lambda) \iff i \in \Delta(\lambda')$. 

\smallskip

We discuss two cases depending on the parity of $\lambda_{i_1+1}- \lambda_{i_1+2}$. 

\smallskip

\noindent
\textbullet \;  {$\lambda_{i_1+1}- \lambda_{i_1+2} \in 2\N$}.\\
Then $i_1 \in \Delta_{\rm bad}(\lambda)$. 
There are two sub-cases depending on whether $i_1 +2$ 
is in $\Delta(\lambda)$ or not:

\begin{itemize}
\item[\textasteriskcentered]  {$i_1 +2 \in \Delta(\lambda)$}. \\
Then, $i_1+2 \in \Delta_{\rm bad}(\lambda)$ (since $\lambda_{i_1+1}- \lambda_{i_1+2} \in 2\N$) 
and $i_1 +2 \in \Delta(\lambda')$. 
Once again, there are two sub-cases:

\begin{itemize} 
\item[1)] {$i_1 +2 \not\in \Delta_{\rm bad}(\lambda')$}.\\ 
Then $| \Delta_{\rm bad}(\lambda')| = | \Delta_{\rm bad}(\lambda)| -2$. 
Moreover,   $(i_1,i_1+2)$ is a good 2-cluster of $\lambda$. 
Indeed,  $ i_1 +2 \not\in \Delta_{\rm bad}(\lambda')$ implies that 
$\lambda_{i_1+3} - \lambda_{i_1+4} \not\in 2\N$.  
On the other hand, $\lambda_{i_1-1} - \lambda_{i_1} =1 \not\in 2\N$ 
(if $i_1=1$ the first condition in Definition \ref{d:cluster} should be omitted).  
But $(i_1,i_1+2)$ is not a 2-cluster of $\lambda'$ since $i_1 \not \in\Delta(\lambda')$.
Hence, we have $|\Sigma(\lambda')|=|\Sigma(\lambda)|- 1$ by Lemma \ref{l:cluster}.

\item[2)]  {$i_1 +2 \in \Delta_{\rm bad}(\lambda')$}.\\
Then $| \Delta_{\rm bad}(\lambda')| = | \Delta_{\rm bad}(\lambda)| -1$. 
Moreover, the unique 2-cluster of $\lambda$ 
which is possibly not a 2-cluster of $\lambda'$ is 
$(i_1,i_1+2)$ but it is bad in this case. 
Indeed,  $\lambda_{i_1+3} - \lambda_{i_1+4} \in 2\N$ 
since $i_1 +2 \in \Delta_{\rm bad}(\lambda')$ 
(and $\lambda'_{i_1+1} - \lambda'_{i_1+2} \not\in 2\N$).
Hence, $|\Sigma(\lambda')|=|\Sigma(\lambda)|$. 
\end{itemize}

\item[\textasteriskcentered] {$i_1 +2 \not\in \Delta(\lambda)$}.\\
Then  $| \Delta_{\rm bad}(\lambda')| = | \Delta_{\rm bad}(\lambda)| -1$. 
Moreover, since $i_1 +2 \not \in \Delta(\lambda)$, then neither $i_1$ 
nor $i_1+2$ belongs to a 2-cluster for $\lambda$. 
Hence $|\Sigma(\lambda)|=|\Sigma(\lambda')|$. 
\end{itemize}

\noindent
\textbullet \; {$\lambda_{i_1+1} - \lambda_{i_1+2} \not\in 2\N$}.\\
In this case, $i_1 \not\in \Delta_{\rm bad}(\lambda)$, $i_1+2 \not \in \Delta(\lambda)$ 
and $i_1+2 \not \in \Delta(\lambda')$. 
Hence $| \Delta_{\rm bad}(\lambda')| = | \Delta_{\rm bad}(\lambda)| $. 
Moreover, neither $i_1$ nor $i_1+2$ belongs to any 2-cluster. 
Hence $|\Sigma(\lambda)|=|\Sigma(\lambda')|$.

\smallskip

In all the cases, we can check with the formula of Definition \ref{d:zPT} 
that $\zPT(\lambda') = \zPT(\lambda) -1$ as desired. 
This concludes the proof of Theorem \ref{t:z}. 
\end{proof}

\section{Counter-examples for \cite[Proposition 3.11]{M08}}  \label{S:err}
From now on, we shall denote 
by $z(\lambda)$ the integer $z_\Mo(\lambda) = \zPT(\lambda)$ 
for $\lambda \in \P_\eps(N)$. 
If $\l$ is a Levi subalgebra of $\g$ and $\O'$ is a rigid nilpotent 
orbit of $\l$, we denote by ${\rm Ind}_{\l}^{\g}(\O')$ the induced nilpotent 
orbit of $\g$ from $\O'$ in $\l$. 

Proposition 3.11 of \cite{M08} asserts that 
if a nilpotent element $e$ associated with the partition 
$\lambda \in \P_\eps(N)$ is induced form a nilpotent 
orbit in a Levi subalgebra $\l$, then $z(\lambda)$ 
is equal to the dimension of the center of $\l$. 
This result is actually incorrect. 
If it were true, it would imply that all 
the sheets containing $e$ share the same dimension (see \cite[Remark 3.12]{M08}). 
But this is wrong. 
Below are some counter-examples 
(see also \cite[Remark 4]{PT}):

\begin{ex}  \label{ex1}
Assume that $\g = \mathfrak{so}(8)$ 
and consider the nilpotent element $e$ of $\g$ 
with partition $\lambda = (3,3,1,1) \in \P_1(8) \smallsetminus \P_1^* (8)$. 
The algorithm yields  $z(\lambda)=2$. 

On the other hand, $e$ is induced from two different ways: 
from the zero orbit in a Levi subalgebra $\l_1$ {\em of type} $(3,1 ; 0)$, 
that is $\l_1 \simeq \mathfrak{gl}_3 \times  \mathfrak{gl}_1 \times 0 $  
(see the definition after \cite[Lemma 3.2]{M08} for the meaning 
of {\em type}),  
and from the zero orbit in a Levi subalgebra $\l_2$ of type $(2 ; 4)$, 
that is $\l_1 \simeq \mathfrak{gl}_2 \times  \mathfrak{so}_4$. 
The first one, $\l_1$, has a center of dimension 2, while 
the second one, $\l_2$, has a center of dimension 1. 
The nilpotent orbit of $e$ has dimension 18 and 
$e$ lies in two different sheets: one of dimension 
$\dim \mathfrak{z}(\l_1) + \dim {\rm Ind}_{\l_1}^{\g}(0)=20$ and one of dimension 
$\dim \mathfrak{z}(\l_2) + \dim {\rm Ind}_{\l_2}^{\g}(0)=19$ 
(here $\mathfrak{z}(\l_i)$ denotes the center of $\l_i$ for $i=1,2$). 
This contradicts Proposition 3.11 of \cite{M08}, 
and also Remark 3.12 of the same paper. 
\end{ex}

\begin{ex}  \label{ex2}
We give now a counter-example in $\mathfrak{sp}(14)$. 
Consider the partition $\lambda=(4,4,2,2,1,1)$ of $\P_{-1}(14)$. 
Here, the algorithm yields $z(\lambda)=2$. 

The corresponding nilpotent element is induced 
from the zero orbit in $\l_1 \simeq \mathfrak{gl}_1 \times \mathfrak{gl}_3 \times \mathfrak{sp}(6)$, 
and from the ridid nilpotent orbit $0 \times \O'$ in $\l_1 \simeq \mathfrak{gl}_2 \times \mathfrak{sp}(10)$ 
where $\O'$ corresponds to the partition $(2,2,2,2,1,1) \in \P_{-1}^*(10)$.  
Again the dimensions of the centers of $\l_1$ and $\l_2$ 
lead to different dimensions, 2 and 1 respectively. 
\end{ex}

The origin of the error can be pined down in the proof 
of \cite[Proposition 3.11]{M08}. 
Let us briefly explain this. 
Until the end of the section, we are in the notations of \cite{M08}. 

At the end of this proof, the assertion 
{\em "Consequently the smallest integer such that one of the situations \textbf{(a)} or \textbf{(b)} of \textbf{Step 1} 
happens in ${\bf d}^{(p)}$ is equal to $i_p$"} is incorrect 
(here ${\bf d}$ is an element of $\P_\eps(N)$). 
And so, the main induction argument of the proof fails. 
We can see that is incorrect in general on an explicit example. 
Consider the partition ${\bf d} = (4,4,3,3,1,1)$ of $\P_{1}(16)$. 
Then the corresponding nilpotent orbit 
is induced from the zero orbit in $\l \simeq \mathfrak{gl}_3 \times \mathfrak{gl}_5 \times 0$ 
and from the rigid nilpotent orbit with partition $(2,2,1,1,1,1)$ 
in $\l\simeq \mathfrak{gl}(4) \times \mathfrak{so}(8)$. 
Consider the second induction. 
In the notations of the proof, we have: 
$S=1$, $i_1=4$, ${\bf d}^{(0)} = {\bf f} = (2,2,1,1,1,1)$, 
${\bf d} = {\bf d}^{(1)} = \widetilde{{\bf d}^{(0)}}$ 
(see \cite[Proposition 3.7]{M08} for the tilda notation). 
Then the smallest integer such that one of the situations \textbf{(a)} or \textbf{(b)} 
of \textbf{Step 1} happens for ${\bf d} ={\bf d}^{(1)}$ is $3 \not = i_1$.

\section{Conclusion}   \label{S:concl}
To summarize, we list below all corrections which have to be taken into account in \cite{M08} 
(the numbering of \cite{M08} is used):

\begin{itemize}
\item Proposition 3.11 (its proof and its statement) is incorrect. 
\item As a consequence Remark 3.12, the sentence "The results of this section specify that, in the classical case, the dimension of a sheet containing a given nilpotent orbit does not depend on the choice of a sheet containing it" 
in \S1.2, and the sentence "Surprisingly, in the classical case, we will notice that if ${\rm Ind}_{\l_1} (\O_{\l_1} ) 
= {\rm Ind}_{\l_2}(\O_{\l_2})$, then $\dim \mathfrak{z}_\g(\l_1) = \dim \mathfrak{z}_\g(\l_2)$" in Remark 2.15, 
are also incorrect. 
\item The proof of Theorem 3.13 is incorrect, since it uses Proposition 3.11. 
Nevertheless, its statement remains valid. 
In particular, Tables 3, 4 and 5 are still correct. 

\noindent
{\em Remark.} 
There are some misprints in Table 5: 
line $2m=48$, the 
partitions are $[7,1^5], [5,3,2^2], [4^2,3,1]$ 
and not $[4^3], [4^2,3,1]$.  
%($[4^3]$ is not an element of $\mathcal{P}_{1}(12)$!)
%\item I take the opportunity to also mention some misprints in this article. 
\end{itemize}


\begin{thebibliography}{DK00} 
\bibitem[BK79]{BoK} W. Borho and H. Kraft,
	{\em \"{U}ber {B}ahnen und deren {D}eformationen bei linearen
	{A}ktionen reduktiver {G}ruppen},
	Comment. Math. Helvetici, {\bf 54} (1979), 61--104. 
\bibitem[B81]{Bor} W. Borho, 
	{\em \"{U}ber {S}chichten halbeinfacher {L}ie-{A}lgebren}, Inventiones Mathematicae,
	{\bf 65} (1981/82), p. 283--317.
\bibitem[CM]{CM} D. Collingwood and W. M. McGovern, 
	{\em Nilpotent orbits in semisimple {L}ie algebras}, 
	Van Nostrand Reinhold Co. New York, {\bf 65} (1993).
\bibitem[K83]{Ke} G. Kempken, 
	{\em Induced conjugacy classes in classical Lie-algebras}, 
	Abh. Math. Sem. Univ. Hamburg, vol. {\bf 53} (1983), 53--83.
\bibitem[M08]{M08} A. Moreau, 
	{\em On the dimension of the sheets of a reductive Lie algebra}, 
	J. Lie Theory, {\bf 18} (2008), n$^\circ$3, 671--696.	
\bibitem[PT12]{PT} A. Premet and L. Topley, 
	{\em Derived subalgebras of centralisers and finite $W$-algebras}, preprint 
	arxiv.org/abs/1301.4653. 
\bibitem[S82]{Sp} N. Spaltenstein, 
	{\em Classes unipotentes et sous-groupes de {B}orel}, 
	Springer-Verlag, Berlin (1982).
\end{thebibliography}
\end{document}